\documentclass[12pt,a4paper]{article}
\usepackage[utf8]{inputenc}
\usepackage{amssymb,amsmath,theorem}

\theorempreskipamount=\medskipamount
\theorempostskipamount=\medskipamount
\newcommand\proof[1]{\textit{Proof#1\pin}}

\newtheorem{theorem}{Theorem}
\newtheorem{lemma}[theorem]{Lemma}
\newtheorem{corollary}[theorem]{Corollary}

\let\eps=\varepsilon
\newcommand\pin{\kern.0833em}
\newcommand{\E}{\mathsf{E}\kern 0.07em}
\newcommand\qed{\ifhmode\unskip\nobreak\fi\quad\ifmmode\square\else\hbox{$\square$}\fi}
\newcommand{\plus}{{\scriptscriptstyle\! +}}

\begin{document}

\begin{center}
{\Large\bfseries Convexity and concavity of $\boldsymbol f$-potentials\\[3pt]
                 (Kolmogorov means)}

\bigskip

{\large\itshape V.\,I. Bakhtin\/$^1$, N.\,A. Tsarev\/$^2$}

\bigskip

$^1$Belarusian State University, Minsk, Belarus (bakhtin@tut.by)

\smallskip

$^2$\pin Belarusian State University, Minsk, Belarus (colya.tsarew2015@yandex.by)
\end{center}

\vspace{-27pt}

\renewcommand{\abstractname}{}
\begin{abstract} \noindent
In the paper we prove criteria for convexity and concavity of $f$-potentials (Kolmogorov means, weighted quasi-arithmetic means), which particular cases are the arithmetic, geometric, harmonic means, the thermodynamic potential (exponential mean), and the $\displaystyle L^{p}$-norm. Then we compute in quadratures all functions $f$ satisfying these criteria.
\end{abstract}


\smallskip


\section{Motivation and results}

Let $(\Omega,\mathfrak A,P)$ be a non-degenerate probability space (i.\,e., one that is a union of two disjoint measurable parts with positive probabilities), and let $f\colon I\to\mathbb R$ be a continuous and strictly monotonic function defined on an interval $I\subset\mathbb R$ (which is possibly infinite in one or both directions).

Denote by $\displaystyle L^{\!c}(\Omega,I)$ the collection of all measurable functions $\varphi\colon \Omega \to I$ with ranges contained compactly in $I$. Obviously, the set $\displaystyle L^{\!c}(\Omega,I)$ is convex.

An $f$-\emph{potential} (or $f$-\emph{mean}) is a functional on $\displaystyle L^{\!c}(\Omega,I)$ defined as
%
\begin{equation} \label{,,1}
 \lambda_f(\varphi) \pin=\pin f^{-1}\bigl(\pin\E\{f(\varphi)\}\bigr) \pin= \pin
 f^{-1}\biggl(\int_\Omega f(\varphi(\omega))\,dP(\omega)\biggr),
 \qquad \varphi\in L^{\!c}(\Omega,I).
\end{equation}
The fact that the range of $\varphi$ is compactly contained in $I$ guarantees convergence of the integral in \eqref{,,1}.

In the case of a discrete probability distribution $P$ (when $\Omega$ is finite) the functional \eqref{,,1} is often called a \emph{weighted quasi-arithmetic mean} \cite{Jarczyk,Pales}; and if, in addition, probabilities of all elements in $\Omega$ are equal, then the term \emph{Kolmogorov---Nagumo---de~Finetti mean} is used (owing to the pioneer papers \cite{Kolmogorov,Nagumo,Finetti}).

Let us give some examples of functionals \eqref{,,1} for different functions $f$.

In the case $f(x) =ax +b$, where $a\ne 0$, we have $\lambda_f(\varphi) = \E\varphi$, the expectation.

If $f(x) =x^p$, where $p>1$, then $\lambda_f(\varphi)$ turns into the $\displaystyle L^{\!\pin p}$-norm $(\pin\E\{\varphi^{\pin p}\})^{1/p}$ (for positive functions $\varphi$).


For $f(x) =1/x$ formula \eqref{,,1} gives the harmonic mean of $\varphi$ (assuming that $\varphi>0$).

For $f(x) =\ln x$ it gives the geometric mean of $\varphi$.

The exponential mean
\begin{equation*}
 \lambda(\varphi) \pin=\pin \ln \E\{e^\varphi\} \pin=\pin
 \ln\biggl(\int_\Omega e^{\varphi(\omega)}\,dP(\omega)\biggr)
\end{equation*}
plays a central role in statistical physics. This is due to the fact that the derivative $\lambda'(\varphi)$ is the Gibbs distribution on $\Omega$ corresponding to the potential function $\varphi$:
\begin{equation*}
 \lambda'(\varphi)[\psi] \pin=\pin \frac{d\lambda(\varphi+t\psi)}{dt}\bigg|_{t=0} =\pin
 \frac{\int e^{\varphi(\omega)}\psi(\omega)\,dP(\omega)}{\int e^{\varphi(\omega)}\,dP(\omega)}.
\end{equation*}
In other words, $\lambda(\varphi)$ is an antiderivative functional for the family of Gibbs distributions. It is often called the thermodynamic potential (like in field theory, where antiderivatives of vector fields are called the potentials).

In the case of an arbitrary function $f$, the derivative $\lambda'_f(\varphi)$ is also a finite (but not probability) measure on $\Omega$ with the density $\rho(\omega) = f'(\lambda_f(\varphi))^{-1} f'(\varphi(\omega))$, because
\begin{equation*}
  \lambda'_f(\varphi)[\psi] \pin=\pin \frac{d\lambda_f(\varphi+t\psi)}{dt}\bigg|_{t=0} =\pin
 \frac{1}{f'(\lambda_f(\varphi))} \int_\Omega f'(\varphi(\omega)) \pin\psi(\omega)\,dP(\omega).
\end{equation*}

In probability theory, $\Gamma(t) =\ln \E\{e^{t\varphi}\}$ is called the cumulant generating function for the random variable $\varphi$. It is used to estimate probabilities of large deviations of the sample mean values of $\varphi$ from the expectation $\E\varphi$.

All of the above convinces of the necessity of studying $f$-means. Below, we call them $f$-potentials (by analogy with statistical physics).

\medskip

It is easy to see that if the function $f(x)$ is replaced with $Af(x)+B$, where $A\ne 0$, then the $f$-potential does not change. Furthermore, it possesses the monotonicity property: the inequality $\psi\ge\varphi$ implies $\lambda_f(\psi)\ge \lambda_f(\varphi)$.

Then, it is natural to try to determine the direction of convexity of the $f$-potential (if any). The purpose of this paper is to prove criteria for convexity and concavity of $f$-potentials given in the following theorems \ref{..1}--\ref{..4}.

\begin{theorem} \label{..1}
Suppose a function\/ $f\colon I\to\mathbb R$ is continuous and increasing. Then convexity of the \emph{$f$-}potential\/ \eqref{,,1} is equivalent to the set of conditions

\smallskip

a) $f$ is convex,

\smallskip

b) $f$ is twice continuously differentiable,

\smallskip

c) the function\/ $h\colon I\to (0,+\infty]$ defined as
\begin{equation} \label{,,2}
 h(x) = \frac{f'(x)}{f''(x)}, \qquad x\in I,
\end{equation}
is concave.
\end{theorem}

\begin{theorem} \label{..2}
Suppose a function\/ $f\colon I\to\mathbb R$ is continuous and increasing. Then convexity of the \emph{$f$-}potential\/ \eqref{,,1} is equivalent to the set of conditions

\smallskip

a) $f$ is convex,

\smallskip

b) $f$ is twice continuously differentiable,

\smallskip

c) the function\/ $H\colon f(I)\to (0,+\infty]$ defined by the implicit formula
\begin{equation} \label{,,3}
 H(f(x)) = \frac{f'(x)^2}{f''(x)}, \qquad x\in I,
\end{equation}
is concave.
\end{theorem}

The only difference between Theorems \ref{..1} and \ref{..2} is in the formulation of conditions c). Below (in Lemma \ref{..6}) it will be proved that conditions c) in these theorems are equivalent to each other, which implies equivalence of Theorems \ref{..1}, \ref{..2}. In addition, it will be established (in Lemma \ref{..5}) that under the conditions of these theorems either $f''(x)>0$ for all $x\in I$ (and then $h$, $H$ are positive and everywhere finite), or $f''(x)\equiv 0$ on all $I$ (and then $f$ is affine, and $h$, $H$ are identically equal to $+\infty$).

Generally speaking, Theorem \ref{..1} looks simpler than Theorem \ref{..2}, because \eqref{,,2} defines the function $h$ explicitly, whereas \eqref{,,3} defines the function $H$ implicitly. Oddly enough, we were unable to find a direct proof of Theorem \ref{..1}, and it will be obtained as a corollary of Theorem \ref{..2}.

Further, Theorems \ref{..1}, \ref{..2} can be easily modified for the cases of a decreasing function $f$ and/or a concave potential $\lambda_f$. In particular, the following theorem is true.

\begin{theorem} \label{..3}
Theorem \ref{..1} stays true under the following combinations of monotonicity and convexity conditions for the functions involved\/$:$

a) if\/ $f$ increases, then the potential\/ $\lambda_f$ is convex\/ $\,\Leftrightarrow\pin$ $f$ is convex and\/ $h$ is concave;

b) if\/ $f$ decreases, then the potential\/ $\lambda_f$ is convex\/ $\,\Leftrightarrow\pin$ $f$ is concave and\/ $h$ is concave;

c) if\/ $f$ decreases, then the potential\/ $\lambda_f$ is concave\/ $\,\Leftrightarrow\pin$ $f$ is convex and\/ $h$ is convex;

d) if\/ $f$ increases, then the potential\/ $\lambda_f$ is concave\/ $\,\Leftrightarrow\pin$ $f$ is concave and\/ $h$ is convex.

\noindent
In other words, under the condition\/ $f\in C^2(I)$, the convexity\/ $($concavity\/$)$ of the $f$-potential is equivalent to the positivity and concavity\/ $($respectively, negativity and convexity\/$)$ of the function\/~$h$.
\end{theorem}

If a function $f$ defined on an interval $I\subset \mathbb R$ is continuous and strictly monotonic, then the inverse function $g = f^{-1}$ is defined on the interval $f(I)$ and is continuous and monotonic in the same sense. It turns out that there is a strict correspondence between the convexity (or concavity) of the potentials $\lambda_f$ and $\lambda_g$.

\begin{theorem} \label{..4}
Suppose $f$ is a continuous and strictly monotonic function defined on the interval\/ $I\subset \mathbb R$, and\/ $g =f^{-1}$ is the inverse function defined on $f(I)$. Then the convexity\/ $($concavity\/$)$ of the functional\/ $\lambda_f$ implies the convexity or concavity\/ of\/ $\lambda_g$. In this case, under the conditions of Theorem \ref{..3}, either both\/ $\lambda_f$ and\/ $\lambda_g$ belong to type b), or both to type c), or one of them is of type a) and the other is of type d).
\end{theorem}

As an illustration (we omit the details of calculations), the table below lists some elementary functions $f$ whose $f$-potentials turn out to be convex or concave. Arrows in the table denote the monotonicity of functions, and arcs denote their convexity or concavity.

In particular, from the fifth and seventh (for $p=-1$) rows of the table it is clear that the geometric mean and the harmonic mean are concave functionals, and from the first row it follows that the exponential mean is a convex functional.


For all other `school' functions ($\sinh$, $\mathop{\mathrm{arsinh}}$, $\tanh$, $\mathop{\mathrm{artanh}}$, and all trigonometric ones not included in the table), restricted to any interval in their domain of definition, the corresponding $f$-potentials are neither convex nor concave.

Solving differential equation \eqref{,,2}, one can compute $f(x)$:
\begin{equation} \label{,,4}
 f(x)\pin=\pin A\int_{x_0}^x \exp\biggl\{\int_{x_0}^s\frac{du}{h(u)}\biggr\}\,ds+B, \quad\ A\ne 0.
\end{equation}
This formula gives an explicit expression for all (except affine) functions $f(x)$ such that the associated potential $\lambda_f$ is convex or concave (respectively, when the function $h(x)$ is positive and concave or negative and convex).

\vspace{-3pt}

$$
\begin{array}{@{\rule{0pt}{15pt}}|c|c|c|c|c|c|c|}
  \cline{1-7}
  \text{No} & f(x) & I & \text{type of}\ f & h(x) & \text{type of}\ h & \text{type of}\ \lambda_f\!\pin \\[3pt]
  \cline{1-7}
  1 & e^x & (-\infty,\infty) & \text{\,a)}\ \nearrow\ \smile\pin & 1 & +\,\frown & \boldsymbol\smile \\
  2 & x^p,\ p>1 & (0,\infty) & \text{\,a)}\ \nearrow\ \smile\pin & \frac{x}{p-1} & +\,\frown & \boldsymbol\smile \\
  3 & \cosh x & (0,\infty) & \text{\,a)}\ \nearrow\ \smile\pin & {}\tanh x & +\,\frown & \boldsymbol\smile \\
  4 & \sec x & (0,\pi/2) & \text{\,a)}\ \nearrow\ \smile\pin & \frac{\sin 2x}{3-\cos 2x} & +\,\frown & \boldsymbol\smile \\[6pt]
  \cline{1-7}
  5 & \ln x & (0,\infty) & \text{\,d)}\ \nearrow\ \frown\pin & -x & -\,\smile & \boldsymbol\frown \\
  6 & \pin x^p,\, p{\,\in\pin}(0,\!\pin 1)\!\pin & (0,\infty) & \text{\,d)}\ \nearrow\ \frown\pin & \frac{x}{p-1} & -\,\smile & \boldsymbol\frown \\
  7 & x^p,\ p<0 & (0,\infty) & \text{\,c)}\ \searrow\ \smile\pin & \frac{x}{p-1} & -\,\smile & \boldsymbol\frown \\
  8 & \mathop{\mathrm{arcosh}}x & (1,\infty) & \text{\,d)}\ \nearrow\ \frown\pin & x^{-1} -x & -\,\smile & \boldsymbol\frown \\
  9 & \coth x & (0,\infty) & \text{\,c)}\ \searrow\ \smile\pin & -\frac{1}{2}\tanh x & -\,\smile & \boldsymbol\frown \\
  10 & \mathop{\mathrm{arcoth}}x & (1,\infty) & \text{\,c)}\ \searrow\ \smile\pin & \frac{1}{2}(x^{-1} -x) \!\pin & -\,\smile & \boldsymbol\frown \\
  11 & \mathop{\mathrm{arcsec}}x & (1,\infty) & \text{\,d)}\ \nearrow\ \frown\pin & -\frac{x^3 -x}{2x^2 -1} & -\,\smile & \boldsymbol\frown \\
  12 & \csc x & (0,\pi/2) & \text{\,c)}\ \searrow\ \smile\pin & -\frac{\sin 2x}{3+\cos 2x} & -\, \smile & \boldsymbol\frown \\
  13 & \mathop{\mathrm{arccsc}}x & (1,\infty) & \text{\,c)}\ \searrow\ \smile\pin & -\frac{x^3 -x}{2x^2 -1} & -\,\smile & \boldsymbol\frown \\[6pt]
  \cline{1-7}
\end{array}
$$

\bigskip

Under the assumptions that $f$ is twice continuously differentiable and the probability distribution $P$ is discrete, Theorem \ref{..1} was proved in \cite[Theorem 2.2]{Pales}, and Theorems \ref{..2}, \ref{..4} were proved in \cite[Theorem 1, Theorem 5]{Jarczyk}. Also under the assumption $f\in C^2(I)$, Theorem \ref{..2} was independently (but later) proved by the second of the authors (see \cite{Tsarev}). This paper provides rigorous proofs for Theorems \ref{..1} -- \ref{..4}, and most of it is devoted to  verifying that the convexity of the $f$-potential implies automatically that $f \in C^2(I)$.


In conclusion, we note that Theorems \ref{..1}--\ref{..4} still stay valid even in the case when the probability measure $P$ on the space $(\Omega,\mathfrak A)$ is only finitely additive, because its countable additivity is never used in the proofs.

\section{Derivation of Theorems \ref{..1}, \ref{..3}, \ref{..4} from Theorem \ref{..2}}

\begin{lemma} \label{..5}
Suppose a function $f\in C^2(I)$ is increasing and convex on the interval\/ $I\subset\mathbb R$, and at least one of the functions\/ $h(x)$, $H(y)$, defined by\/ \eqref{,,2}, \eqref{,,3}, is concave. Then there is an alternative: either $f''(x)\ne 0$ on all of\/ $I$, or $f''(x) \equiv 0$.
\end{lemma}

\proof. Under the conditions of lemma the functions $h$ and $H$ take values in $(0,+\infty]$. The one of them that is concave has a convex subgraph. Therefore, it is either everywhere finite or identically equal to $+\infty$. Respectively, either $f''(x)\ne 0$ on all $I$, or $f''(x)\equiv 0$. \qed

\begin{lemma} \label{..6}
Suppose that a function $f\in C^2(I)$ is strictly monotonic, and $f''(x)\ne 0$ on all\/ $I$. Then in the case of increasing $f$ the convexity\/ $($concavity\/$)$ of the function\/ $h(x)$ from\/ \eqref{,,2} is equivalent to the convexity\/ $($concavity\/$)$ of the function\/ $H(y)$ from\/ \eqref{,,3}$;$ and in the case of decreasing $f$ the convexity\/ $($concavity\/$)$ of\/ $h(x)$ is respectively equivalent to the concavity\/ $($convexity\/$)$ of\/ $H(y)$.
\end{lemma}

\proof. From \eqref{,,2}, \eqref{,,3} it follows that
\begin{equation} \label{,,5}
  f''(x) =\frac{f'(x)}{h(x)}, \qquad f''(x) = \frac{f'(x)^2}{H(f(x))}.
\end{equation}
If one of the functions $h(x)$ or $H(y)$ is convex or concave, then its right-hand derivative ($h'_\plus(x)$ or $H'_\plus(y)$, respectively) is well defined on all its domain. In both cases, it follows from \eqref{,,5} that the function $f''(x)$ has the right-hand derivative $f'''_\plus(x)$ on the interval $I$, and then from \eqref{,,2}, \eqref{,,3} we obtain the right-hand derivatives
\begin{gather*}
  h'_\plus(x) \pin=\pin 1 -\frac{f'(x)f'''_\plus(x)}{f''(x)^2}, \\[6pt]
  H'_\plus(y)\Bigr|_{y=f(x)} \pin=\pin \biggl[2f'(x) -\frac{f'(x)^2f'''_\plus(x)}{f''(x)^2}\biggr] \frac{1}{f'(x)} \pin=\pin 2 -\frac{f'(x)f'''_\plus(x)}{f''(x)^2}.
\end{gather*}
Subtracting them yields the identity
\begin{equation*}
  H'_\plus(f(x)) -h'_\plus(x) \pin\equiv\pin 1, \quad\  x\in I.
\end{equation*}
It implies that in the case of increasing $f$ the functions $H'_\plus(y)$ and $h'_\plus(x)$ have the same monotonicity, and so $H(y)$ and $h(x)$ have the same convexity (if at least one of them possesses one). Similarly, in the case of decreasing $f$ the directions of monotonicity and convexity of these functions turn out to be opposite. \qed

\medskip

From Lemmas \ref{..5} and \ref{..6} it is clear that conditions c) in Theorems \ref{..1} and \ref{..2} are equivalent to each other. Therefore, these theorems are equivalent, and it is sufficient to prove only one of them (namely, Theorem \ref{..2}).

\medskip

Let us now assume that Theorem \ref{..1} is true, and derive from it Theorems \ref{..3}, \ref{..4}.

Obviously, item a) of theorem \ref{..3} duplicates Theorem \ref{..1}. Items b), c), d) are derived from a) by means of replacing an increasing convex function $f(x)$ with a decreasing concave function $g(x) =-f(x)$, with a decreasing convex function $g(x) =f(-x)$, and with an increasing concave function $g(x) =-f(-x)$, respectively. As a result of these replacements, the potential $\lambda_g$ and the corresponding function $h_g(x) =g'(x)/g''(x)$ (similar to \eqref{,,2}) turn out to be as follows:

b) the potential $\lambda_g(\varphi) =\lambda_f(\varphi)$ is convex, and the function $h_g(x) =h(x)$ is concave;

c) the potential $\lambda_g(\varphi) =-\lambda_f(-\varphi)$ is concave, and $h_g(x) = -h(-x)$ is convex;

d) the potential $\lambda_g(\varphi) =-\lambda_f(-\varphi)$ is concave, and $h_g(x) = -h(-x)$ is convex.

Thus, theorem \ref{..3} is proved.

\medskip

Now, let us prove Theorem \ref{..4} under the conditions of Theorem \ref{..1} (or, equivalently, part a) of Theorem \ref{..3}), namely, when the function $f(x)$ is increasing and convex, and $h(x)$ is concave. In this setting the inverse to $f(x)$ function $g(y)$ is increasing and concave. Differentiate it twice:
\begin{gather*}
 g'(y) =\frac{1}{f'(x)} =\frac{1}{f'(g(y))}, \\[3pt]
 g''(y) =-\frac{1}{f'(g(y))^2}\pin f''(g(y))\pin \frac{1}{f'(g(y))} =-\frac{f''(g(y))}{f'(g(y))^3},
\end{gather*}
and, using these derivatives, express the function $H(y)$ from \eqref{,,3}:
\begin{equation} \label{,,6}
 H(y) = \frac{f'(g(y))^2}{f'(g(y))} = -\frac{g'(y)}{g''(y)}.
\end{equation}

By Lemma \ref{..5}, either $f''(x) \equiv 0$ or $f''(x) >0$ on the entire interval $I$.

In the first case the functions $f(x)$ and $g(y)$ are affine, and the corresponding potentials $\lambda_f(\varphi) =\lambda_g(\varphi) =\E\{\varphi\}$ are linear.

If $f''(x) >0$, then by Lemma \ref{..6} the function $H(y)$, as well as $h(x)$, is concave, and from \eqref{,,6} it follows that the function $h_g(y) =g'(y)/g''(y)$ is convex. Thus, $g$ satisfies the conditions of part d) of Theorem \ref{..3}, and hence the potential $\lambda_g$ is concave.

Under the conditions of parts b), c), d) of Theorem \ref{..3}, Theorem \ref{..4} is proved in the same way.

\section{Proof of the sufficient part of Theorem \ref{..2}}

Let a function $f\in C^2(I)$ be increasing and convex, and the function $h(x)$ from \eqref{,,2} be concave. Then by Lemma \ref{..5} either $f''(x) \equiv 0$ or $f''(x) >0$ on the whole of $I$.

In the first case the function $f(x)$ is affine, the potential $\lambda_f(\varphi) = \E\{\varphi\}$ is linear and hence is convex.

To prove the convexity of the functional $\lambda_f$ in the case when $f''(x) >0$, it is sufficient to verify that
\begin{equation*}
 \frac{d^{\pin 2}\lambda_f(\varphi+t\psi)}{dt^2}\biggr|_{t=0} \pin\ge\pin 0 \quad \text{for all}\ \ \varphi \in L^{\!c}(\Omega,I),\ \, \psi\in L^\infty(\Omega).
\end{equation*}

By \eqref{,,1} we have
\begin{equation*}
  \lambda_f(\varphi +t\psi) \pin=\pin f^{-1}\bigl(\pin\E\{f(\varphi+t\psi)\}\bigr).
\end{equation*}
Therefore,
\begin{gather} \notag
 \frac{d\lambda_f(\varphi +t\psi)}{dt} \pin=\pin
 \frac{\E\{f'(\varphi+t\psi)\pin\psi\}}{f'(\lambda_f(\varphi +t\psi))}, \\[9pt] \notag
 \frac{d^{\pin 2}\lambda_f(\varphi +t\psi)}{dt^2}\biggr|_{t=0} =\pin
 \frac{\E\{f''(\varphi)\pin\psi^2\}}{f'(\lambda_f(\varphi))} -
 \frac{f''(\lambda_f(\varphi))}{f'(\lambda_f(\varphi))^2} \pin
 \frac{\E\{f'(\varphi)\pin\psi\}^2}{f'(\lambda_f(\varphi))} \pin=\pin \\[6pt] \notag
 \pin=\pin\frac{f''(\lambda_f(\varphi))}{f'(\lambda_f(\varphi))^3}
 \biggl(\!\pin \frac{f'(\lambda_f(\varphi))^2}{f''(\lambda_f(\varphi))}
 \, \E\bigl\{f''(\varphi)\pin\psi^2\bigr\} - \E\{f'(\varphi)\pin\psi\}^2\biggr) \pin=\pin \\[6pt]
 \pin=\pin \frac{f''(\lambda_f(\varphi))}{f'(\lambda_f(\varphi))^3}
 \biggl(\!\pin H\bigl(\E\{f(\varphi)\}\bigr)\pin \E\bigl\{f''(\varphi)\pin\psi^2\bigr\} -
 \E\biggl\{\!\pin \frac{f'(\varphi)}{f''(\varphi)^{1/2}} \cdot f''(\varphi)^{1/2}\psi
 \biggr\}^{\! 2}\,\biggr),  \label{,,7}
\end{gather}
where the function $H(y)$ is taken from \eqref{,,3}. By Lemma \ref{..6} $H(y)$ is concave, and, by virtue of Jensen's inequality,
\begin{equation} \label{,,8}
  H\bigl(\E\{f(\varphi)\}\bigr) \ge \E\bigl\{H(f(\varphi))\bigr\} =
  \E\biggl\{\frac{f'(\varphi)^2}{f''(\varphi)}\biggr\}.
\end{equation}
Replacing in \eqref{,,7} the expression $H\bigl(\E\{f(\varphi)\}\bigr)$ by the right-hand side of \eqref{,,8} and applying the  Cauchy--Bunyakovsky--Schwarz inequality, we see that the expression \eqref{,,7} is nonnegative, and hence the functional $\lambda_f$ is convex.

\section{Proof of the necessary part of Theorem \ref{..2}}

The necessary part of Theorem \ref{..2} is much more difficult. Its proof takes up the rest of the paper.

Once and for all we fix a partition of the space $\Omega$ into two measurable parts with positive probabilities:
\begin{equation} \label{,,9}
 \Omega = \Omega_0 \sqcup \Omega_1, \ \ \text{where}\ \, p_0 =P(\Omega_0) >0 \ \, \text{and}\ \,
 p_1 =P(\Omega_1) >0.
\end{equation}
In what follows we use only functions $\displaystyle \varphi \in L^{\!c}(\Omega,I)$ taking constant values on $\Omega_0$ and $\Omega_1$.

\begin{lemma} \label{..7}
If the \emph{$f$-}potential\/ \eqref{,,1} is convex, then the function\/ $f$ is convex.
\end{lemma}

\proof{} is by contradiction. Assume that the function $f$ is not convex. Then there exist points $x_0,x_1\in I$ and a number $\theta\in (0,1)$ such that
\begin{equation*}
 f\big((1-\theta)x_0 +\theta x_1\big) > (1-\theta)f(x_0) +\theta f(x_1).
\end{equation*}

Consider the function
\begin{equation} \label{,,10}
 g(x) = f(x)  -kx, \quad \text{where}\ \ k =\frac{f(x_1) -f(x_0)}{x_1-x_0}.
\end{equation}
It satisfies the relations
\begin{equation*}
 g(x_0) =g(x_1) < g\big((1-\theta)x_0 +\theta x_1\big).
\end{equation*}
Denote by $x_*$ the leftmost maximum point of $g(x)$ on $[x_0,x_1]$. Then
\begin{equation} \label{,,11}
 \begin{cases}
 g(x) <g(x_*), & \text{if}\ \ x_0\le x < x_*, \\[2pt]
 g(x)\le g(x_*), & \text{if}\ \ x_*< x \le x_1.
 \end{cases}
\end{equation}

Using \eqref{,,9}, let us define a one-parameter family of functions
\begin{equation*}
 \varphi_t(\omega) =
  \begin{cases}
    x_* +tp_1, & \text{if}\ \ \omega\in \Omega_0, \\[2pt]
    x_* -tp_0, & \text{if}\ \ \omega\in \Omega_1.
  \end{cases}
\end{equation*}
Evidently, $\displaystyle \varphi_t\in L^{\!c}(\Omega,I)$ for small $t$. And if $t\ne 0$, then \eqref{,,1}, \eqref{,,10}, \eqref{,,11} imply that
\begin{align*}
 \lambda_f(\varphi_t) &=
 f^{-1}\big(p_0f(x_* +tp_1) +p_1 f(x_* -tp_0)\big) \\[3pt]
 &= f^{-1}\big(p_0g(x_* +tp_1) +p_0k(x_* +tp_1) + p_1g(x_* -tp_0) + p_1k(x_* -tp_0)\big) \\[3pt]
 &< f^{-1}\big(p_0g(x_*) +p_0kx_* +p_1g(x_*) +p_1kx_*\big) = f^{-1}\big(g(x_*) +kx_*\big) =x_*,
\end{align*}
whereas $\lambda_f(\varphi_0) =x_*$. Hence $\lambda_f(\varphi_t)$ has a strict local maximum at the point $t=0$, which contradicts the convexity of $\lambda_f$. \qed

\medskip

The convexity of $f(x)$ implies the existence of the right-hand derivative
\begin{equation*}
 f'_\plus(x) =\lim_{t\to 0+} \frac{f(x+t) -f(x)}{t}, \qquad x \in I,
\end{equation*}
which is nondecreasing and right continuous. And the monotonicity of $f'_\plus(x)$ implies that for almost all $x\in I$ there exists the ordinary derivative
\begin{equation} \label{,,12}
 \frac{df'_\plus(x)}{dx} \pin=\pin \lim_{t\to 0} \frac{f'_\plus(x+t) -f'_\plus(x)}{t}.
\end{equation}

Let $I_0$ be the set of all points $x\in I$ where the derivative \eqref{,,12} exists. Then its complement $I\setminus I_0$ has zero Lebesgue measure. Let $J =f(I)$ and $J_0 =f(I_0)$. Since the right-hand derivative $f'_\plus(x)$ is monotonic, it is locally bounded. It follows that $f(x)$ is locally Lipschitz continuous, and hence the set $J\setminus J_0 = f(I\setminus I_0)$ has zero measure.

\begin{lemma} \label{..8}
If the \emph{$f$-}potential\/ \eqref{,,1} is convex, then
\begin{equation} \label{,,13}
 \limsup_{t\to 0} \frac{|f'_\plus(x+t) -f'_\plus(x)|}{|t|} < \infty \quad \text{for each}\ \
 x\in I.
\end{equation}
\end{lemma}

\begin{corollary} \label{..9}
If the \emph{$f$-}potential\/ \eqref{,,1} is convex, then\/ $f\in C^1(I)$.
\end{corollary}

\proof. Fix an $x\in I$. Consider all pairs of numbers $y_0,y_1\in J$ satisfying the equality
\begin{equation} \label{,,14}
  f(x) =p_0y_0 +p_1y_1,
\end{equation}
where $p_0$ and $p_1$ are taken from \eqref{,,9}. It is clear from \eqref{,,14} that $y_0$ and $y_1$ are affinely dependent on each other. So we can choose $y_0,y_1\in J_0$ (in fact, those $y_0\in J_0$ for which $y_1\notin J_0$ have zero measure).

Set $x_0 =f^{-1}(y_0)$ and $x_1 =f^{-1}(y_1)$. Then \eqref{,,14} takes the form
\begin{equation} \label{,,15}
  f(x) =p_0f(x_0) +p_1f(x_1), \quad x_0,x_1\in I_0.
\end{equation}

Consider a function $\displaystyle \varphi\in L^{\!c}(\Omega, I)$ defined as
\begin{equation*}
 \varphi(\omega) =
  \begin{cases}
    x_0, & \text{if}\ \  \omega\in \Omega_0, \\[2pt]
    x_1, & \text{if}\ \  \omega\in \Omega_1.
  \end{cases}
\end{equation*}
Let $g(t) =\lambda_f(\varphi +t)$. The function $g(t)$ is defined in a neighborhood of zero, is convex and increasing. Therefore, it has a nondecreasing and positive right-hand derivative
$g'_\plus(t)$.

The right-hand derivative of the composition $f(g(t))$ is equal to $f'_\plus(g(t))\pin g'_\plus(t)$. On the other hand, by \eqref{,,1} we have
\begin{equation*}
 f(g(t)) =\E\{f(\varphi +t)\} =p_0f(x_0+t) +p_1f(x_1+t),
\end{equation*}
which implies
\begin{equation} \label{,,16}
 f'_\plus(g(t))\pin g'_\plus(t) \pin=\pin p_0f'_\plus(x_0+t) +p_1f'_\plus(x_1+t).
\end{equation}

\medskip

Since the functions $f'_\plus(g(t))$ and $g'_\plus(t)$ are positive and nondecreasing,
\begin{equation} \label{,,17}
  \bigl| f'_\plus(g(t))\pin g'_\plus(0) -f'_\plus(g(0))\pin g'_\plus(0)\bigr| \pin\le\pin
 \bigl| f'_\plus(g(t))\pin g'_\plus(t) -f'_\plus(g(0))\pin g'_\plus(0)\bigr|.
\end{equation}
Substitution of \eqref{,,16} into the right-hand side of \eqref{,,17} yields the estimate
\begin{equation} \label{,,18}
 \bigl| f'_\plus(g(t))\pin g'_\plus(0) -f'_\plus(g(0))\pin g'_\plus(0)\bigr| \pin\le\pin
 \bigl|p_0 (f'_\plus(x_0+t) -f'_\plus(x_0)) +p_1 (f'_\plus(x_1+t) -f'_\plus(x_1))\bigr|.
\end{equation}
Combining \eqref{,,18} with differentiability of the function $f'_\plus(x)$ at the points $x_0,x_1\in I_0$, we see that
\begin{equation} \label{,,19}
 \limsup_{t\to 0} \frac{|f'_\plus(g(t)) -f'_\plus(g(0))|}{|t|} < \infty.
\end{equation}

\medskip

Recall that $g(t)$ is convex and increasing. Hence, the difference $g(t)-g(0)$ has the same order of smallness as $t$. Besides, by \eqref{,,1} and \eqref{,,15},
\begin{equation*}
 g(0) =\lambda_f(\varphi) = f^{-1}\bigl(\pin\E\{f(\varphi)\}\bigr) =
 f^{-1}\bigl(p_0f(x_0) +p_1f(x_1)\bigr) = f^{-1}(f(x)) =x.
\end{equation*}
Therefore \eqref{,,19} implies \eqref{,,13}. \qed

\medskip

Thus, the convexity of $\lambda_f$ implies the convexity and continuous differentiability of $f$. In this situation the derivative $f'(x) =f'_\plus(x)$ is nondecreasing and differentiable on the subset of full measure $I_0\subset I$. Consequently, there exists a second derivative $f''(x)$ at each point $x\in I_0$, and we can define
\begin{equation} \label{,,20}
 H(f(x)) = \frac{f'(x)^2}{f''(x)}, \qquad x\in I_0.
\end{equation}
The function $H(y)$ is defined on the set $J_0 = f(I_0)$, which has full measure in $J =f(I)$.

\begin{lemma} \label{..10}
If the \emph{$f$-}potential\/ \eqref{,,1} is convex, and points\/ $y_0,\,y_1$ and\/ $y=p_0y_0 +p_1y_1$ $($where the numbers\/ $p_0,\,p_1$ are taken from\/ \eqref{,,9}$)$ lie in\/ $J_0$, then
\begin{equation} \label{,,21}
 H(y) \ge p_0\pin H(y_0) +p_1H(y_1).
\end{equation}
\end{lemma}

\proof. If $H(y) =+\infty$ then \eqref{,,21} is trivial. So we further assume that $H(y)
<+\infty$. Let $x_0,\,x_1,\,x$ denote the $f$-preimages of $y_0,\,y_1,\,y$, respectively. Under the conditions of Lemma \ref{..10} they all lie in $I_0$ and satisfy the equality
\begin{equation} \label{,,22}
 f(x) =p_0f(x_0) +p_1f(x_1).
\end{equation}
Moreover, the finiteness of $H(y)$ and the convexity of $f$ imply that $f''(x) >0$.

Consider the function
\begin{equation*}
 \varphi(\omega) =
  \begin{cases}
    x_0, & \text{if}\ \ \omega\in \Omega_0, \\[2pt]
    x_1, & \text{if}\ \ \omega\in \Omega_1.
  \end{cases}
\end{equation*}

 \medskip\noindent
From \eqref{,,1} and \eqref{,,22} it follows that $\lambda_f(\varphi) =x$.

Take an arbitrary function $\psi$ with constant values $\psi_0$ and $\psi_1$ on the sets
$\Omega_0$ and $\Omega_1$ respectively, and calculate the derivatives
\begin{gather} \notag
 \frac{d\lambda_f(\varphi +t\psi)}{dt} \pin=\pin
 \frac{p_0f'(x_0 +t\psi_0)\pin\psi_0 +p_1f'(x_1+t\psi_1)\pin\psi_1}{f'(\lambda_f(\varphi +t\psi))},
 \\[6pt] \label{,,23}
 \frac{d^{\pin 2}\lambda_f(\varphi+t\psi)}{dt^2}\biggr|_{t=0} =\pin
 \frac{p_0f''(x_0)\pin \psi_0^2 +p_1f''(x_1)\pin \psi_1^2}{f'(x)} -
 \frac{f''(x)[\pin p_0f'(x_0)\pin\psi_0 +p_1f'(x_1)\pin\psi_1]^2}{f'(x)^3}.
\end{gather}

Due to the convexity of $\lambda_f$ the right-hand side of \eqref{,,23} is nonnegative. Recall that the derivatives $f'(x)$, $f'(x_0)$, $f'(x_1)$, $f''(x)$ are strictly positive, and $f''(x_0)$, $f''(x_1)$ are nonnegative. If $f''(x_0) =0$, then taking $\psi_0 =1$ and $\psi_1 =0$ we get a negative right-hand side \eqref{,,23}, which cannot be. Hence $f''(x_0) >0$. Similarly, $f''(x_1) >0$.

Now let us take
\begin{equation} \label{,,24}
 \psi_0 =\frac{f'(x_0)}{f''(x_0)}, \qquad \psi_1 =\frac{f'(x_1)}{f''(x_1)},
\end{equation}
and substitute \eqref{,,24} into \eqref{,,23}:
\begin{gather} \notag
 \frac{1}{f'(x)}\biggl[p_0\frac{f'(x_0)^2}{f''(x_0)}+p_1\frac{f'(x_1)^2}{f''(x_1)}\biggr] -
 \frac{f''(x)}{f'(x)^3}\biggl[p_0\frac{f'(x_0)^2}{f''(x_0)}+p_1\frac{f'(x_1)^2}{f''(x_1)}\biggr]^2
 = \\[6pt] \label{,,25}
 =\pin \frac{f''(x)}{f'(x)^3}
 \biggl[p_0\frac{f'(x_0)^2}{f''(x_0)}+p_1\frac{f'(x_1)^2}{f''(x_1)}\biggr]
 \biggl(\frac{f'(x)^2}{f''(x)}-p_0\frac{f'(x_0)^2}{f''(x_0)}-p_1\frac{f'(x_1)^2}{f''(x_1)}\biggr).
\end{gather}
Since the whole expression \eqref{,,25} is nonnegative, its rightmost factor must be nonnegative. On the other hand, this factor coincides with $H(y) -p_0\pin H(y_0) -p_1H(y_1)$. So Lemma \ref{..10} is proved. \qed

\begin{lemma} \label{..11}
If the \emph{$f$-}potential\/ \eqref{,,1} is convex, then either\/ $H(y) =+\infty$ for all\/ $y\in J_0$, or\/ $H(y)<+\infty$ for all\/ $y\in J_0$.
\end{lemma}

\proof. Assume that $H(y_0) =+\infty$ at some point $y_0\in J_0$. Then from \eqref{,,21} it follows that $H(y) =+\infty$ at almost all points $y =p_0y_0 +p_1y_1$, where $y_1\in J$ (except those ones for which $y\notin J_0$ or $y_1\notin J_0$). Successive repetition of this argument extends the equality $H(y) =+\infty$ to some subset of full measure in the whole of $J$.

On the other hand, if there exists a point $y\in J_0$ such that $H(y)<+\infty$, then \eqref{,,21} imply that $H(y_0) <+\infty$ and $H(y_1)<+\infty$ for all $y_0,y_1\in J_0$ satisfying the condition $y =p_0y_0 +p_1y_1$. Obviously, such points $y_0,\,y_1$ fill up a set of full measure in a neighborhood of $y$. This contradicts the conclusion of the previous paragraph. Therefore, the conditions $H(y_0) =+\infty$ and $H(y)<+\infty$ are incompatible for any $y_0,y\in J_0$. \qed

\begin{lemma} \label{..12}
Suppose\/ $J =(a,b)$. Under the conditions of Lemma \ref{..11}, if the function\/ $H(y)$ takes finite values on the whole of\/ $J_0$, then
\begin{equation}\label{,,26}
 \frac{-2H(y_0)}{b-y_0} \pin\le\pin \frac{H(y) -H(y_0)}{y-y_0} \pin\le\pin \frac{2H(y_0)}{y_0-a} \quad\ \text{for all}\ \ y_0,y\in J_0.
\end{equation}
\end{lemma}

\begin{corollary} \label{..13}
If the function\/ $H(y)$ takes finite values on\/ $J_0$, then its restriction to\/ $J_0$ is locally bounded and locally Lipschitz continuous.
\end{corollary}


\proof. From \eqref{,,21} it follows that if the segment $[y_0,y_1]$ is divided by a point $y$ in proportion $p_0\!:\! p_1$ or $p_1\!:\! p_0$, and all the three points lie in $J_0$, then
\begin{equation} \label{,,27}
 \frac{H(y)-H(y_0)}{y-y_0} \pin\ge\pin \frac{H(y_1)-H(y_0)}{y_1-y_0} \pin\ge\pin \frac{H(y_1)-H(y)}{y_1-y}.
\end{equation}
Since $p_0+p_1=1$, without loss of generality we assume $p_1\ge 1/2$.

First let us consider the case $y\in (y_0,b)$. Define a sequence of points
\begin{equation*}
 y_n =y_0 + (y-y_0)p_1^{-n+1}, \qquad n\in \mathbb N.
\end{equation*}
It is easy to see that $y_n$ divides the segment $[y_0,y_{n+1}]$ in proportion $p_1\!:\! p_0$.

In the case $b<+\infty$ choose a natural number $n$ such that
\begin{equation*}
 \frac{b-y_0}{2} \le y_n -y_0 < b-y_0.
\end{equation*}
Obviously, the points $y_1,\,\dots,\,y_n$ belong to $J_0$ for almost every $y\in J_0\cap (y_0,b)$.
Iteratively applying the left inequality \eqref{,,27}, we get
\begin{equation} \label{,,28}
 \frac{H(y)-H(y_0)}{y-y_0} \pin=\pin \frac{H(y_1)-H(y_0)}{y_1-y_0} \pin\ge\pin \ \cdots \
 \pin\ge\pin \frac{H(y_n)-H(y_0)}{y_n-y_0} \pin\ge\pin \frac{-2H(y_0)}{b-y_0}.
\end{equation}
This implies the validity of the left inequality \eqref{,,26} for almost all $y\in (y_0,b)$. If $b=+\infty$, then the chain of inequalities \eqref{,,28} becomes infinite, and yields the same left inequality \eqref{,,26} in the limit.

Similarly, continue the chain of inequalities \eqref{,,28} to the left using the points
\begin{equation*}
 y_{-1} =y+ (y_0-y)p_1^{-1}, \qquad y_{-n} =y_0+ (y_{-1}-y_0)p_1^{-n+1}, \ \ n\in \mathbb N.
\end{equation*}
For almost every $y\in J_0\cap (y_0,b)$ all these points lie in $J_0$ (except those out of $J$) and, when $y$ varies, they fill up a set of full measure in the interval $(a,y_0)$. Iteratively applying the right inequality \eqref{,,27} and at the very end \eqref{,,28}, we get
\begin{equation*}
 \cdots\ \pin\ge\pin  \frac{H(y_0)-H(y_{-n})}{y_0-y_{-n}} \pin\ge\pin \ \cdots\ \pin\ge\pin
 \frac{H(y_0)-H(y_{-1})}{y_0-y_{-1}} \pin\ge\pin  \frac{H(y) -H(y_0)}{y -y_0}
 \pin\ge\pin \frac{-2H(y_0)}{b-y_0}.
\end{equation*}
Replacing here each $y_{-n}$ by $y$, we obtain the left inequality \eqref{,,26} for almost all points $y\in J_0\cap (a,y_0)$. And taking into account the previous paragraph, \eqref{,,26} is proved for almost all $y\in J_0$. Now we have to extend it to \emph{all} $y\in J_0$.

Suppose $y_0,y_1\in J_0$. Then, as proved, for almost all $y\in J_0\cap (y_0,y_1)$
\begin{align*}
 \frac{H(y_1) -H(y_0)}{y_1-y_0} &\pin=\pin \frac{H(y) -H(y_0)}{y-y_0} \frac{y-y_0}{y_1-y_0} +
 \frac{H(y) -H(y_1)}{y-y_1} \frac{y_1-y}{y_1-y_0} \pin\ge\pin \\[6pt]
   &\pin\ge\pin
 \frac{-2H(y_0)}{b-y_0}\frac{y-y_0}{y_1-y_0} + \frac{-2H(y_1)}{b-y_1}\frac{y_1-y}{y_1-y_0}
 \ \longrightarrow\ \frac{-2H(y_0)}{b-y_0} \quad \text{as}\ \ y\to y_1
\end{align*}
(where it does not matter which is greater, $y_0$ or $y_1$). Hence, the left inequality \eqref{,,26} is true for all $y_0,y\in J_0$.

The right inequality \eqref{,,26} is proved in the same way. \qed

\begin{lemma} \label{..14}
Suppose a function\/ $g\colon I\to\mathbb R$ satisfies the condition
\begin{equation} \label{,,29}
 \limsup_{t\to 0} \frac{|g(x+t) -g(x)|}{|t|} < \infty \quad \text{for all}\ \ x\in I,
\end{equation}
and the limit\/ \eqref{,,29} is equal to zero for almost all\/ $x\in I$. Then\/ $g \equiv \mathrm{const}$.
\end{lemma}

\proof. It is sufficient to consider a finite interval $I$. Let $I_0$ denote the set of those points $x\in I$ at which the limit \eqref{,,29} is equal to zero, and let $I_n$ denote the set of those points $x\in I$ at which the limit \eqref{,,29} belongs to $(0,n)$. By the condition, the set $I_0$ has full measure in $I$, whereas for each $n>0$ the set $I_n$ has measure zero.

Fix an $\eps>0$. For every $x\in I_0$, using \eqref{,,29}, choose a neighborhood $O(x)$ such that the oscillation of $g$ over $O(x)$ does not exceed $\eps |O(x)|$. Then for each $n>0$ define an open set $U_n\supset I_n$ of Lebesgue measure less than $\eps/2^n$, and for every $x\in I_n$, using \eqref{,,29}, choose a neighborhood $O(x)\subset U_n$ such that the oscillation of $g$ over it does not exceed $2n |O(x)|$.

The neighborhoods $O(x)$ constructed in the previous paragraph form a covering of the interval $I$. One can extract from it a countable subcovering that splits into two parts, say $\mathcal A$ and $\mathcal B$, each consisting of pairwise disjoint neighborhoods (by means of Besicovitch covering theorem \cite{Besicovitch}).

Let us estimate the sum of oscillations of $g$ over the neighborhoods $O(x)\in \mathcal A$:
\begin{equation} \label{,,30}
 \sum_{O(x)\in \mathcal A} \mathop{\mathrm{osc}}(g,O(x)) \pin\le\pin
 \sum_{n=0}^\infty \sum_{\substack{O(x)\in \mathcal A,\\ x\in I_n}} \mathop{\mathrm{osc}} (g,O(x)) \pin\le\pin \eps|I| +\sum_{n=1}^\infty \frac{2n\eps}{2^n}.
\end{equation}
The sum of oscillations of $g$ over the neighborhoods $O(x)\in \mathcal B$ can be estimated in exactly the same way. It follows that the oscillation of $g$ over the interval $I$ does not exceed twice the right-hand side of \eqref{,,30}, and due to the arbitrariness of $\varepsilon$ it is zero. \qed

\begin{lemma} \label{..15}
If the \emph{$f$-}potential\/ \eqref{,,1} is convex, then\/ $f\in C^2(I)$.
\end{lemma}

\proof. By Lemma \ref{..11} there are two alternatives: either $f''(x) =0$ for all $x\in I_0$, or $f''(x) >0$ for all $x\in I_0$.

Assume that the first alternative takes place. Then, by virtue of Lemma \ref{..8}, the function $g(x) =f'(x)$ satisfies the conditions of Lemma \ref{..14}, and according to it $f'(x) \equiv \mathrm{const}$. Obviously, in this case $f\in C^2(I)$.

Now assume that the second alternative takes place. Then, by virtue of Corollary \ref{..13}, the function $H(y)$ can be continuously extended to $J$. Accordingly, the function
\begin{equation*}
 f''(x) =\frac{f'(x)^2}{H(f(x))}
\end{equation*}
(see \eqref{,,20}) extends continuously from $I_0$ to $I$. Fix an $x_0\in I$ and set
\begin{equation} \label{,,31}
 g(x) = f'(x) -f'(x_0) -\int_{x_0}^x \frac{f'(t)^2}{H(f(t))}\, dt.
\end{equation}
This function $g(x)$ satisfies the conditions of Lemma \ref{..14}, and according to it $g(x) \equiv 0$. Then \eqref{,,31} shows that $f\in C^2(I)$. \qed

\begin{lemma} \label{..16}
If the \emph{$f$-}potential\/ \eqref{,,1} is convex, then\/ $H(y)$ from\/ \eqref{,,3} is concave.
\end{lemma}

\proof. Since $f\in C^2(I)$, the set $I_0$ coincides with $I$, and the set $J_0$ coincides with $J$. In this setting Lemma \ref{..11} states that either $H(y) \equiv +\infty$ or $H(y) <
+\infty$ on the entire $J$. In the first case Lemma \ref{..16} is trivial. In the second case the function $H(y)$ is continuous on $J$, and by Lemma \ref{..10},
\begin{equation} \label{,,32}
  H(p_0y_0 +p_1y_1) \ge p_0\pin H(y_0) +p_1H(y_1) \quad \text{for all}\ \ y_0,y_1\in J.
\end{equation}
But this is not sufficient, because the coefficients $p_0,\pin p_1$ are \emph{fixed} here.

Consider the function
\begin{equation*}
  g(y) =H(y) -H(y_0) -\frac{H(y_1)-H(y_0)}{y_1-y_0}(y-y_0).
\end{equation*}
Obviously, $g(y_0) =g(y_1) =0$, and \eqref{,,32} imply
\begin{equation*}
 g(p_0y_0 +p_1y_1)\ge 0.
\end{equation*}
Repeatedly dividing segments in proportion $p_1\!:\! p_0$, we can construct an everywhere dense subset $Y\subset [y_0,y_1]$ such that the function $g(y)$ is nonnegative on $Y$. By continuity this function is nonnegative on the entire $[y_0,y_1]$. This is equivalent to \eqref{,,32} with arbitrary coefficients $p_0,\, p_1>0$ that add up to one. Hence $H(y)$ is concave. \qed

\medskip

Now the necessary part of Theorem \ref{..2} follows from Lemmas \ref{..7}, \ref{..15}, \ref{..16}.



\end{document}